\documentclass[12pt,twoside,reqno]{amsart}
\usepackage[colorlinks=true,citecolor=blue]{hyperref}
\usepackage{mathptmx, amsmath, amssymb, amsfonts, amsthm, mathptmx, enumerate, color,mathrsfs}
\usepackage{graphicx}

\usepackage{multirow}
\usepackage{epstopdf}
\usepackage{multicol}
\usepackage{algorithm}
\usepackage{algorithmic}
\usepackage{epstopdf}

\newcommand{\zero}{\mathbf{0}}
\newcommand{\R}{\mathbb{R}}
\newcommand{\co}{\mathrm{co}}

\newcommand{\N}{\mathbb{N}} 

\newcommand{\C}{\mathbb{C}}
\newcommand{\Q}{\mathbb{Q}}
\newcommand{\PP}{\mathbb{P}}

\newtheorem{theorem}{Theorem}[section]

\theoremstyle{definition}

\newtheorem{remark}[theorem]{Remark}
\numberwithin{equation}{section}

\begin{document}
\setcounter{page}{1}

\vspace*{2.0cm}
\title[Optimality conditions, KKT and neural networks]
{KKT-based optimality conditions for neural network approximation}
\author[V. Peiris, N. Sukhorukova, J. Ugon]{ Vinesha Peiris$^{1}$, Nadezda Sukhorukova$^{2*}$, Julien Ugon$^3$}
\maketitle
\vspace*{-0.6cm}

\begin{center}
{\footnotesize

$^1$Centre for Optimisation and Decision Science, Curtin University, Australia\\
$^2$Department of Mathematics, Swinburne  University of Technology, Australia\\
$^3$School of Information Technology, Deakin  University, Australia

}\end{center}

\vskip 4mm {\footnotesize \noindent {\bf Abstract.}
In this paper, we obtain necessary optimality conditions for neural network approximation. We consider neural networks in Manhattan ($l_1$ norm) and Chebyshev ($\max$ norm). The optimality conditions are based on neural networks with at most one hidden layer. We reformulate nonsmooth unconstrained optimisation problems as larger dimension constrained problems with smooth objective functions and constraints. Then we use KKT conditions to develop the necessary conditions and present the optimality conditions in terms of convex analysis and convex sets.

 \noindent {\bf Keywords.}
min-max approximation, KKT, optimality conditions, neural networks

 \noindent {\bf 2020 Mathematics Subject Classification.}
41A50, 26A27, 65D10. }

\renewcommand{\thefootnote}{}
\footnotetext{ $^*$Corresponding author.
\par
E-mail address: m.peiris@curtin.edu.au (V. Peiris), nsukhorukova@swin.edu.au (N. Sukhorukova), julien.ugon@deakin.edu.au (J. Ugon).
\par

}

\section{Introduction}

The theory of Artificial Neural Networks (ANNs) or just neural networks appeared as an extension of work by A.~Kolmogorov and V.~Arnold~\cite{Arnold57,Kol57}, who were working on the 13th problem of Hilbert. Their celebrated representation theorem and further developments lead to an efficient adaptation to function approximation, established in~\cite{Cybenko,Hornik1991,LeshnoPinkus1993,pinkus_1999}. A detailed overview of optimisation approaches for studying neural networks can be found in~\cite{SunOptDeepLearning}. 

In this paper, we approach neural networks purely from the point of view of optimisation. Therefore, we are minimising the amount of ``purely'' neural network vocabulary. Essentially, neural networks provide accurate approximations to (continuous) functions. In most practical problems, function values are defined at a finite number of points (discretisation points). 

Our approach is to formulate the corresponding approximation problems as optimisation problems. These problems are nonsmooth and nonconvex. We reformulate these nonsmooth problems as equivalent smooth ones by introducing extra variables and constraints (standard approach) and then formulate the corresponding KKT conditions, that are necessary optimality conditions, but, in general, they are not sufficient.

The rest of the paper is organised as follows. In section~\ref{sec:preliminary}, we present the necessary knowledge from the theory of neural networks and formulate the corresponding optimisation problems. Then in section~\ref{sec:approximation}, we formulate the corresponding KKT conditions and interpret them from the point of view of convex analysis. Finally, in section~\ref{sec:conclusions}, we present the conclusions and possible further research directions.

\section{Preliminaries}\label{sec:preliminary}
\subsection{Neural networks and approximation}
Neural networks can be seen as approximation prototypes. In this paper, we focus on the simplest neural networks (with at most one hidden layer). Without too many details, the structure of the corresponding approximations is as follows:
\begin{enumerate}
\item single hidden layer (one hidden layer):
\begin{equation}\label{eq:singlelayer}
\sum_{j=1}^{n}a_j\sigma((w_1^j,w_2^j,\dots,w_d^j){\bf T_i}+w_0^j),   
\end{equation}
\item no hidden layer (zero hidden layers, index $j$ is redundant):
\begin{equation}\label{eq:nolayer}
\sigma((w_1,w_2,\dots,w_d){\bf T_i}+w_0).
\end{equation}
\end{enumerate}
In these formulations,  ${\bf T_i\in \R^d}$, $i=1,\dots,N$ are the discretisation points and $$a_1,\dots,a_n,w_0^j, w_1^j\dots, w_d^j,~j=1,\dots,n$$
are the approximation parameters ($w_0^j, w_1^j\dots, w_d^j$, $j=1,\dots,n$ are also known as network weights), $\sigma:\R\to\R$ is the activation function (univariate). 

\begin{remark}
    In the case of no hidden layers~(\ref{eq:nolayer}), we assume for simplicity that $a_1=1$. This coefficient is ``absorbed'' by the activation function without changing the smoothness.  Since in~(\ref{eq:nolayer}) there is only one set of weights, we omit the upper index for simplicity.
\end{remark}

The approximation parameters are the decision variable of the corresponding optimisation problems. It is clear that the number of decision variables is significantly higher in the case of one hidden layer.

The activation function is predefined and not subject to optimisation. A common choice for activation functions are sigmoidal functions, ReLU and Leaky ReLU~\cite{Goodfellow2016}. In~\cite{Hornik1991} and~\cite{pinkus_1999} the authors showed that any nonpolynomial activation function can be used for accurate approximation constructions. 

The goal is to minimise a loss function, which is simply an error function. The approximation errors (deviations) are calculated at the discretisation points and then combined in the chosen loss function. Common choices of loss functions are:
\begin{enumerate}
\item least squares (sum of squares of deviation at discretisation points, $L_2^2$);
\item manhattan loss function (sum of absolute deviations of errors, $L_1$);
\item uniform loss function (maximum of absolute deviations of errors, $L_{\max}$).
\end{enumerate}
Least squares is a popular choice: in the case of smooth activation functions, the corresponding optimisation problems are smooth. In this paper, however, we focus on nonsmooth loss functions, namely, $L_1$ and $L_{\max}$. We also assume that the activation function is smooth, leaving ReLU and Leaky ReLU for future research directions. In the case of nonsmooth activation functions, one has to apply nonsmooth optimisation techniques~\cite{DemyanovDixon,minimax}, since the gradients of the functions appearing in the KKT conditions are not defined due to their nonsmoothness. 

\begin{remark}
In~\cite{PeirisSukhNonsmoothNN} the authors obtained equivalent necessary optimality conditions for the uniform-based loss function (no hidden layer) using nonsmooth optimisation techniques. In the current paper, we show that these conditions can be obtained using the KKT conditions (without applying nonsmooth optimisation techniques) when the activation function is smooth. In addition, in this paper, we develop the optimality conditions for $L_1$. 
\end{remark}

\subsection{Uniform-based loss function}\label{ssec:uniform}

The optimisation problem is as follows:
\begin{equation}\label{eq:free_disc}
 {\text{minimise}}  \max_{\substack{{\bf T}_i\\  i=1,\dots,N}}  \left| \sum_{j=1}^{n}a_j\sigma({\bf w}^j{\bf T}_i+w_{0}^j)- f({\bf T}_i)\right|, \end{equation}
 \begin{equation}\label{eq:free_disc_c0nstraints}
  {\text{subject to}}~ X\in \R^{n+nd+n},
\end{equation}
where 
\begin{itemize}
    \item $N$ is the number of discretisation points,
    \item ${\bf T}_i\in Q$, $i=1,\dots,N$ are $d$-dimensional vectors (dataset points) where $Q$ is a $d$-dimensional hypercube,
    \item $n$ is the maximal number of components in the sum (that is, the number of nodes in the hidden layer), $n=1$ in the case of the absence of the hidden layers (index~$j$ is redundant) and $n>1$ for the case with one hidden layer,
    \item vectors~${\bf w}^j=(w_1^j,\dots,w_d^j)$ are the weights where $w_k^j$, $j=1,\dots,n$, $k=0,\dots,d$, and coefficients~$a_j$, $j=1,\dots,n$ are the decision variables. All the decision variables are placed in a single vector, 
    $$X=(a_1, \dots, a_n, w_{1}^1, \dots, w_{1}^n,\dots, w_{d}^{1},\dots, w_{d}^{n}, w_{0}^1, \dots, w_{0}^n),$$
    \item $\sigma$ is the activation function and, 
    \item $f$ is the original function (function to approximate).
\end{itemize}

An equivalent formulation of~(\ref{eq:free_disc})-(\ref{eq:free_disc_c0nstraints}), where all the objective function and the constraint functions are smooth, is as follows:
\begin{equation}\label{eq:free_disc_uniform}
{\text{minimise}}  ~z 
\end{equation}
subject to
 \begin{gather}\label{eq:free_disc_c0nstraints_uniform}
 \sum_{j=1}^{n}a_j\sigma({\bf w}^j{\bf T}_i+w_{0}^j)- f({\bf T}_i)\leq z\\
 f({\bf T}_i)-\sum_{j=1}^{n}a_j\sigma({\bf w}^j{\bf T}_i+w_{0}^j)\leq z\\  
  X\in \R^{n+nd+n},~i=1,\dots,N.
\end{gather}
The objective function and the constraint functions are smooth and therefore KKT conditions can be used to formulate the  necessary optimality conditions. One more decision variable $z$ was introduced. The total number of constraints is $2N$, which is large in most practical problems.

\begin{remark}
    In the case of no hidden layers and monotonous activation function the problem is quasiconvex and can be efficiently solved using, for example, a bisection method~\cite{PeirisRoshchinsSukhorukova2024}.
\end{remark}

\subsection{Manhattan-based loss function}
The optimisation problem is as follows:
\begin{equation}\label{eq:free_discl1}
 {\text{minimise}}  ~\sum_{i=1}^{N} \left| \sum_{j=1}^{n}a_j\sigma({\bf w}^j{\bf T}_i+w_{0}^j)- f({\bf T}_i)\right|, \end{equation}
 \begin{equation}\label{eq:free_disc_c0nstraintsl1}
  {\text{subject to}}~ X\in \R^{n+nd+n},
\end{equation}
where the components are the same as in section~\ref{ssec:uniform}.

An equivalent formulation of~(\ref{eq:free_discl1})-(\ref{eq:free_disc_c0nstraintsl1}), where all the objective function and the constraint functions are smooth, is as follows:
\begin{equation}\label{eq:free_disc_l11}
{\text{minimise}}  ~\sum_{i=1}^{N}z_i 
\end{equation}
subject to
 \begin{gather}\label{eq:free_disc_c0nstraints_l11}
 \sum_{j=1}^{n}a_j\sigma({\bf w}^j{\bf T}_i+w_{0}^j)- f({\bf T}_i)\leq z_i\\
 f({\bf T}_i)-\sum_{j=1}^{n}a_j\sigma({\bf w}^j{\bf T}_i+w_{0}^j)\leq z_i\\  
  X\in \R^{n+nd+n},~i=1,\dots,N.
\end{gather}
Similar to the uniform-based loss function case, the objective function and the constraint functions are smooth and therefore KKT conditions can be used to formulate necessary optimality conditions. Extra $N$ decision variable $z_i$, $i=1,\dots,N$ were introduced. The total number of constraints is $2N$.

In the next section, we consider the simplest cases of neural networks: no hidden layer and one single layer. At first glance, this is an obvious limitation. However, it is important to start with simpler cases. Moreover, due to~\cite{Cybenko}, even one hidden layer lead to accurate approximations. The only requirement is to choose nonpolynomial functions as activation functions~\cite{LeshnoPinkus1993, pinkus_1999}. At the same time, the number of hidden nodes may be large~\cite{negres}.

\section{Neural networks approximation}\label{sec:approximation}

In this section, we re-numerate the indices of the discretisation points. This is done in order to express the optimality conditions in terms of positive and negative deviation. As we will see shortly, in the case of uniform-based  approximation, the optimality conditions are based on maximal absolute deviation points only, while in the case of Manhattan-based approximation, one has to take into account all the discretisation points, including the point where the deviation is zero.

The renumeration is as follows:
\begin{itemize}
    \item indices $i=1,\dots,n_1$, $i\in \PP$ correspond to positive maximal absolute deviation points, such that  
    $$\sigma({\bf w}^j{\bf T}_i+w_{0}^j)- f({\bf T}_i)>0;$$
    \item indices $i=n_1+1,\dots,n_1+n_2$, $i\in \N$ correspond to negative maximal absolute deviation points, such that  
    $$\sigma({\bf w}^j{\bf T}_i+w_{0}^j)- f({\bf T}_i)<0;$$
    \item remaining indices $i=n_1+n_2+1,\dots,n_1+n_2+n_3$, $i\in \C$ correspond to the points whose deviation is zero.
\end{itemize}
Note that $n_1+n_2+n_3=N$.

\subsection{Uniform-based loss function}
\subsubsection{No hidden layer}
In this case, the decision variables are 
$$(w_0,w_1,\dots,w_d,z).$$ 
Then KKT conditions are as follows:
\begin{equation}\label{eq:KKT1}
\begin{bmatrix}
0 \\
0 \\
\vdots \\
0\\
1
\end{bmatrix}
+\sum_{i=1}^{n_1}\lambda_i
\begin{bmatrix}
\begin{bmatrix}
1 \\
{\bf T}_i \\
\end{bmatrix}\sigma'({\bf w}{\bf T}_i+w_0)\\
-1
\end{bmatrix}
+\sum_{i=n_1+1}^{n_1+n_2}\lambda_i
\begin{bmatrix}
\begin{bmatrix}
-1 \\
-{\bf T}_i \\
\end{bmatrix}\sigma'({\bf w}{\bf T}_i+w_0)\\
-1
\end{bmatrix}
=\zero_{d+2},
\end{equation}
where 
\begin{equation}\label{eq:KKT11}
    \lambda_i\geq 0,~i=1,\dots,n_1+n_2.
\end{equation}
and $\sigma'$ is the derivative of the univariate activation function $\sigma$. 
Note that due to the complementarity conditions, $$\lambda_i=0,~i=n_1+n_2+1,\dots,n_1+n_2+n_3.$$ 
From the last row of~(\ref{eq:KKT1}) one can see that
$$\sum_{i=1}^{n_1+n_2}\lambda_i=1,$$
while from the first row
$$\sum_{i=1}^{n_1}\lambda_i=\sum_{i=n_1+1}^{n_1+n_2}\lambda_i=0.5.$$
Therefore, KKT can be written as follows:
\begin{equation}\label{eq:KKTco1}
    \co\left\{\begin{bmatrix}
1 \\
{\bf T}_i \\
\end{bmatrix}\sigma'({\bf w}{\bf T}_i+w_0),~i\in\PP \right\}\cap\co\left\{\begin{bmatrix}
1 \\
{\bf T}_i \\
\end{bmatrix}\sigma'({\bf w}{\bf T}_i+w_0),~i\in\N \right\}\neq \emptyset.
\end{equation}
The condition~(\ref{eq:KKTco1}) is a necessary optimality condition.

\subsubsection{One hidden layer}
In this case, the decision variables are 
$$(w_0^1,w_1^1,\dots,w_d^1,\dots w_0^n,w_1^n,\dots,w_d^n,a_1,\dots,a_n,z),$$
where $a_1,\dots,a_n$ are the coefficients from the representation~(\ref{eq:singlelayer}).
Then KKT conditions are as follows:
\begin{equation}\label{eq:KKTsingle}
\begin{bmatrix}
0 \\
0 \\
\vdots \\
0\\
1
\end{bmatrix}
+\sum_{i=1}^{n_1}\lambda_i
\begin{bmatrix}
a_1\begin{bmatrix}
1 \\
{\bf T}_i \\
\end{bmatrix}\sigma'({\bf w_1}{\bf T}_i+w_0)\\
\vdots\\
a_n\begin{bmatrix}
1 \\
{\bf T}_i 
\end{bmatrix}\sigma'({\bf w_n}{\bf T}_i+w_0)\\
\begin{bmatrix}
\sigma({\bf w}_1{\bf T}_i+w_1^0)\\
\vdots\\
\sigma({\bf w}_n{\bf T}_i+w_n^0)\\
\end{bmatrix}\\
-1
\end{bmatrix}
+\sum_{i=n_1+1}^{n_1+n_2}\lambda_i
\begin{bmatrix}
a_1\begin{bmatrix}
-1 \\
-{\bf T}_i \\
\end{bmatrix}\sigma'({\bf w_1}{\bf T}_i+w_0)\\
\vdots\\
a_n\begin{bmatrix}
-1 \\
-{\bf T}_i 
\end{bmatrix}\sigma'({\bf w_n}{\bf T}_i+w_0)\\
\begin{bmatrix}
-\sigma({\bf w}_1{\bf T}_i+w_1^0)\\
\vdots\\
-\sigma({\bf w}_n{\bf T}_i+w_n^0)\\
\end{bmatrix}\\
-1
\end{bmatrix}
=\zero_{d+2+n},
\end{equation}
\begin{equation}\label{eq:KKT11single}
    \lambda_i\geq 0,~i=1,\dots,n_1+n_2.
\end{equation}
Similarly to the situation without a hidden layer, using the last coordinate, one can see that
$$\sum_{i=1}^{n_1+n_2}\lambda_i=1.$$
Therefore, a necessary optimality condition can be presented as follows:
\begin{equation}\label{eq:necessary_uniform_single_geom}
\co\left\{\begin{bmatrix}
a_1\begin{bmatrix}
1 \\
{\bf T}_i \\
\end{bmatrix}\sigma'({\bf w_1}{\bf T}_i+w_0)\\
\vdots\\
a_n\begin{bmatrix}
1 \\
{\bf T}_i 
\end{bmatrix}\sigma'({\bf w_n}{\bf T}_i+w_0)\\
\begin{bmatrix}
\sigma({\bf w}_1{\bf T}_i+w_1^0)\\
\vdots\\
\sigma({\bf w}_n{\bf T}_i+w_n^0)\\
\end{bmatrix}
\end{bmatrix},~
i\in\PP
\right\}\cap\co\left\{
\begin{bmatrix}
a_1\begin{bmatrix}
1 \\
{\bf T}_i \\
\end{bmatrix}\sigma'({\bf w_1}{\bf T}_i+w_0)\\
\vdots\\
a_n\begin{bmatrix}
1 \\
{\bf T}_i 
\end{bmatrix}\sigma'({\bf w_n}{\bf T}_i+w_0)\\
\begin{bmatrix}
\sigma({\bf w}_1{\bf T}_i+w_1^0)\\
\vdots\\
\sigma({\bf w}_n{\bf T}_i+w_n^0)\\
\end{bmatrix}
\end{bmatrix},~
i\in\N\right\}\neq\emptyset,
\end{equation}
where set $\PP$ is the set of positive maximal absolute  deviation points, while set $\N$ is the set of negative maximal absolute deviation points. 

Similarly to the case without hidden layers, the necessary optimality conditions are based on maximal absolute deviation points. The situation in the case of one hidden layer is significantly more complex.

\subsection{Manhattan-based loss function}
For Manhattan-based approximation, we use the following notation:
\begin{itemize}
    \item ${\bf E}_k$ is a $k$-th dimensional column vector, whose components are equal to~1;
    \item ${\bf O}_k$ is a $k$-th dimensional column vector, whose components are equal to~0;
   \item  ${\bf J}_{ki}$ is a $k$-th dimensional column vector, whose components are equal to~0, while the $i$-th component is equal to~1.
\end{itemize}
\subsubsection{No hidden layer}

In this case, the decision variables are $$(w_0,w_1,\dots,w_d,z_1,\dots,z_N).$$ 
Then KKT conditions are as follows:
\begin{equation}\label{eq:KKT3}
\begin{bmatrix}
0 \\
0 \\
\vdots \\
0\\
{\bf E}_{n_1}\\
{\bf E}_{n_2}\\
{\bf E}_{n_3}
\end{bmatrix}
+\sum_{i=1}^{n_1}\lambda_i
\begin{bmatrix}
\begin{bmatrix}
1 \\
{\bf T}_i \\
\end{bmatrix}\sigma'({\bf w}{\bf T}_i+w_0)\\
-{\bf J}_{n_1i}\\
 {\bf O}_{n_2} \\
  {\bf O}_{n_3}\\
\end{bmatrix}
+\sum_{i=n_1+1}^{n_1+n_2}\lambda_i
\begin{bmatrix}
\begin{bmatrix}
-1 \\
-{\bf T}_i \\
\end{bmatrix}\sigma'({\bf w}{\bf T}_i+w_0)\\
{\bf O}_{n_1}\\
 -{\bf J}_{n_2 i} \\
  {\bf O}_{n_3}\\
\end{bmatrix}
+{\bf D}
=\zero_{d+1+N},
\end{equation}
where
\begin{equation}\label{eq:additional}
{\bf D}=\sum_{i=n_1+n_2+1}^{N}\lambda_i
\begin{bmatrix}
\begin{bmatrix}
1 \\
{\bf T}_i \\
\end{bmatrix}\sigma'({\bf w}{\bf T}_i+w_0)\\
{\bf O}_{n_1}\\
 {\bf O}_{n_2} \\
  {\bf J}_{n_3 i}\\
\end{bmatrix}
+\sum_{i=1+n_1+n_2+n_3}^{n_1+n_2+2n_3}\lambda_i
\begin{bmatrix}
\begin{bmatrix}
-1 \\
-{\bf T}_i \\
\end{bmatrix}\sigma'({\bf w}{\bf T}_i+w_0)\\
{\bf O}_{n_1}\\
 {\bf O}_{n_2} \\
  {\bf J}_{n_3 i}\\
\end{bmatrix},
\end{equation}
$\sigma'$ is the derivative of the univariate activation function $\sigma$,
\begin{equation}\label{eq:KKT31}
    \lambda_i\geq 0,~i=1,\dots,n_1+n_2+2n_3.
\end{equation}
Therefore, $$\lambda_i=1,~i=1,\dots,n_1+n_2$$
and $$\lambda_{n_1+n_2+i}+\lambda_{n_1+n_2+n_3+i}=1,~i=1,\dots,n_3.$$
Then the necessary optimality conditions are 
\begin{equation}
\sum_{i=1}^{n_1}\begin{bmatrix}
1 \\
{\bf T}_i \\
\end{bmatrix}\sigma'({\bf w}{\bf T}_i+w_0)-\sum_{i=n_1+1}^{n_1+n_2}\begin{bmatrix}
1 \\
{\bf T}_i \\
\end{bmatrix}\sigma'({\bf w}{\bf T}_i+w_0)\in {\bf \Q}
\end{equation}
where
\begin{equation}
\Q=\sum_{i=n_1+n_2+1}^{N} \co
\left\{
\begin{bmatrix}
1 \\
{\bf T}_i \\
\end{bmatrix}\sigma'({\bf w}{\bf T}_i+w_0),-\begin{bmatrix}
1 \\
{\bf T}_i \\
\end{bmatrix}\sigma'({\bf w}{\bf T}_i+w_0)
\right\}
    \end{equation}

In the case when there is no discretisation point, where the deviation is zero, set $\Q=\zero_{d+1}$ and the necessary condition can be simplified:
$$\sum_{i=1}^{n_1}\begin{bmatrix}
1 \\
{\bf T}_i \\
\end{bmatrix}\sigma'({\bf w}{\bf T}_i+w_0)=\sum_{i=n_1+1}^{n_1+n_2}\begin{bmatrix}
1 \\
{\bf T}_i \\
\end{bmatrix}\sigma'({\bf w}{\bf T}_i+w_0).$$
Therefore, in this case the polytopes that characterise the KKT conditions are singletons and hence they are much simpler to verify.

\subsubsection{One hidden layer}

In this case, the decision variables are $$(w_0^1,w_1^1,\dots,w_d^1,\dots w_0^n,w_1^n,\dots,w_d^n,a_1\dots,a_n,z_1,\dots,z_N).$$ 
Then KKT conditions are as follows:
\begin{equation}\label{eq:KKT3hidden}
\begin{bmatrix}
0 \\
0 \\
\vdots \\
0\\
{\bf E}_{n_1}\\
{\bf E}_{n_2}\\
{\bf E}_{n_3}
\end{bmatrix}
+\sum_{i=1}^{n_1}\lambda_i
\begin{bmatrix}
a_1\begin{bmatrix}
1 \\
{\bf T}_i \\
\end{bmatrix}\sigma'({\bf w}^1{\bf T}_i+w_0^1)\\
\vdots\\
a_n\begin{bmatrix}
1 \\
{\bf T}_i \\
\end{bmatrix}\sigma'({\bf w}^n{\bf T}_i+w_0^n)\\
\begin{bmatrix}
\sigma({\bf w}^1{\bf T}_i+w^1_0)\\
\vdots\\
\sigma({\bf w}^n{\bf T}_i+w^n_0)\\
\end{bmatrix}\\
-{\bf J}_{n_1i}\\
 {\bf O}_{n_2} \\
  {\bf O}_{n_3}\\
\end{bmatrix}
+\sum_{i=n_1+1}^{n_1+n_2}\lambda_i
\begin{bmatrix}
a_1\begin{bmatrix}
-1 \\
-{\bf T}_i \\
\end{bmatrix}\sigma'({\bf w}^1{\bf T}_i+w_0^1)\\
\vdots\\
a_n\begin{bmatrix}
-1 \\
-{\bf T}_i \\
\end{bmatrix}\sigma'({\bf w}^n{\bf T}_i+w_0^n)\\
\begin{bmatrix}
-\sigma({\bf w}_1^T{\bf T}_i+w_1^0)\\
\vdots\\
-\sigma({\bf w}_n^T{\bf T}_i+w_n^0)\\
\end{bmatrix}\\
{\bf O}_{n_1}\\
 -{\bf J}_{n_2 i} \\
  {\bf O}_{n_3}\\
\end{bmatrix}
+{\bf D}
=\zero_{d+1+N+n},
\end{equation}
where
\begin{equation}\label{eq:additional}
{\bf D}=\sum_{i=n_1+n_2+1}^{N}\lambda_i
\begin{bmatrix}
a_1\begin{bmatrix}
1 \\
{\bf T}_i \\
\end{bmatrix}\sigma'({\bf w}^1{\bf T}_i+w_0^1)\\
\vdots\\
a_n\begin{bmatrix}
1 \\
{\bf T}_i \\
\end{bmatrix}\sigma'({\bf w}^n{\bf T}_i+w_0^n)\\
\begin{bmatrix}
\sigma({\bf w}^1{\bf T}_i+w_1^0)\\
\vdots\\
\sigma({\bf w}^n{\bf T}_i+w_n^0)\\
\end{bmatrix}\\
{\bf O}_{n_1}\\
 {\bf O}_{n_2} \\
  {\bf J}_{n_3 i}\\
\end{bmatrix}
+\sum_{i=1+n_1+n_2+n_3}^{n_1+n_2+2n_3}\lambda_i
\begin{bmatrix}
a_1\begin{bmatrix}
-1 \\
-{\bf T}_i \\
\end{bmatrix}\sigma'({\bf w}^1{\bf T}_i+w_0^1)\\
\vdots\\
a_n\begin{bmatrix}
-1 \\
-{\bf T}_i \\
\end{bmatrix}\sigma'({\bf w}^n{\bf T}_i+w_0^n)\\
\begin{bmatrix}
-\sigma({\bf w}^1{\bf T}_i+w^1_0)\\
\vdots\\
-\sigma({\bf w}^n{\bf T}_i+w^n_0)\\
\end{bmatrix}\\
{\bf O}_{n_1}\\
 {\bf O}_{n_2} \\
  {\bf J}_{n_3 i}\\
\end{bmatrix}
\end{equation}

\begin{equation}\label{eq:KKT31}
    \lambda_i\geq 0,~i=1,\dots,n_1+n_2+2n_3.
\end{equation}
Therefore, $$\lambda_i=1,~i=1,\dots,n_1+n_2$$
and $$\lambda_{n_1+n_2+i}+\lambda_{n_1+n_2+n_3+i}=1,~i=1,\dots,n_3.$$
Then the necessary optimality conditions are 
\begin{equation}
\sum_{i=1}^{n_1}
\begin{bmatrix}
a_1\begin{bmatrix}
1 \\
{\bf T}_i \\
\end{bmatrix}\sigma'({\bf w}^1{\bf T}_i+w_0^1)\\
\vdots\\
a_n\begin{bmatrix}
1 \\
{\bf T}_i \\
\end{bmatrix}\sigma'({\bf w}^n{\bf T}_i+w_0^n)\\
\begin{bmatrix}
\sigma({\bf w}^1{\bf T}_i+w^1_0)\\
\vdots\\
\sigma({\bf w}^n{\bf T}_i+w^n_0)\\
\end{bmatrix}\\

\end{bmatrix}
+\sum_{i=n_1+1}^{n_1+n_2}
\begin{bmatrix}
a_1\begin{bmatrix}
-1 \\
-{\bf T}_i \\
\end{bmatrix}\sigma'({\bf w}^1{\bf T}_i+w_0^1)\\
\vdots\\
a_n\begin{bmatrix}
-1 \\
-{\bf T}_i \\
\end{bmatrix}\sigma'({\bf w}^n{\bf T}_i+w_0^n)\\
\begin{bmatrix}
-\sigma({\bf w}^1{\bf T}_i+w^1_0)\\
\vdots\\
-\sigma({\bf w}^n{\bf T}_i+w^n_0)\\
\end{bmatrix}\\
\end{bmatrix}
\in {\bf \Q},
\end{equation}
where
\begin{equation}
\Q=\sum_{i=n_1+n_2+1}^{N} 
\begin{bmatrix}
\co
\left\{
\begin{bmatrix}
1 \\
{\bf T}_i \\
\end{bmatrix}a_1\sigma'({\bf w}^1{\bf T}_i+w_0),-\begin{bmatrix}
1 \\
{\bf T}_i 
\end{bmatrix}a_1\sigma'({\bf w}^1{\bf T}_i+w_0)
\right\}\\
\vdots\\
\co
\left\{
\begin{bmatrix}
1 \\
{\bf T}_i \\
\end{bmatrix}a_n\sigma'({\bf w}^n{\bf T}_i+w_0),-\begin{bmatrix}
1 \\
{\bf T}_i \\
\end{bmatrix}a_n\sigma'({\bf w}^n{\bf T}_i+w_0)
\right\}\\
\sigma({\bf w}^1{\bf T}_i+w_0^1)\co\{1,-1\}\\
\vdots\\
\sigma({\bf w}^1{\bf T}_i+w_0^1)\co\{1,-1\}\\
    \end{bmatrix}.
    \end{equation}
Similarly to the case without hidden layers, the conditions can be simplified if none of the discretisation points corresponds to zero deviation:
$$
\sum_{i=1}^{n_1}
\begin{bmatrix}
a_1\begin{bmatrix}
1 \\
{\bf T}_i \\
\end{bmatrix}\sigma'({\bf w}^1{\bf T}_i+w_0^1)\\
\vdots\\
a_n\begin{bmatrix}
1 \\
{\bf T}_i \\
\end{bmatrix}\sigma'({\bf w}^n{\bf T}_i+w_0^n)\\
\begin{bmatrix}
\sigma({\bf w}^1{\bf T}_i+w^1_0)\\
\vdots\\
\sigma({\bf w}^n{\bf T}_i+w^n_0)\\
\end{bmatrix}\\

\end{bmatrix}
=\sum_{i=n_1+1}^{n_1+n_2}
\begin{bmatrix}
a_1\begin{bmatrix}
1 \\
{\bf T}_i \\
\end{bmatrix}\sigma'({\bf w}^1{\bf T}_i+w_0^1)\\
\vdots\\
a_n\begin{bmatrix}
1 \\
{\bf T}_i \\
\end{bmatrix}\sigma'({\bf w}^n{\bf T}_i+w_0^n)\\
\begin{bmatrix}
\sigma({\bf w}^1{\bf T}_i+w^1_0)\\
\vdots\\
\sigma({\bf w}^n{\bf T}_i+w^n_0)\\
\end{bmatrix}\\
\end{bmatrix}.
$$

Therefore, similar to the case without hidden layers (no zero deviation discretisation points), the polytopes that characterise the KKT conditions are singletons.

\section{Conclusions}\label{sec:conclusions}
In this paper, we developed necessary optimality conditions for neural network approximation for at most one hidden layer (uniform and Manhattan loss functions). The nonsmooth optimisation unconstrained problems are formulated as smooth optimisation problems with smooth constraints and therefore the corresponding KKT can be used. 

In the case of the uniform loss function, the conditions are based on the overlapping of polytopes whose vertices correspond to the maximal deviation points. Essentially, only the maximal deviation points contribute to these conditions. 

In the case of the Manhattan loss function, the situation is different. First of all, the optimality conditions, formulated in terms of convex polytopes, include all the discretisation points (not just maximal deviation points). Second, the conditions can be simplified in the case when none of the discretisation points corresponds to the deviation equal to zero: the corresponding polytopes are singletons.    


\end{document}